\documentclass[12pt]{article}

\newtheorem{theorem}{Theorem}[section]
\newtheorem{corollary}[theorem]{Corollary}
\newtheorem{lemma}[theorem]{Lemma}
\newtheorem{proposition}[theorem]{Proposition}
\newtheorem{remark}[theorem]{Remark}

\DeclareSymbolFont{AMSa}{U}{msa}{m}{n}
\DeclareMathSymbol{\boxtimes}{\mathbin}{AMSa}{"02}
\def\fin{\nolinebreak \raise0.2em\hbox{\framebox{}}}
\def\proof{{\noindent\it Proof: }}

\def\dual{ ^{\rm v}}
\def\T{{ ^t}}

\def\LL{{\cal L}}
\def\II{{\cal I}}

\def\quot#1#2{\mbox{\raise 0.1cm\hbox{$#1$}{\big /}\raise -0.1cm\hbox{$#2$}}}

\def\Tor{{\cal T}or }
\def\ext{\mathrm{Ext}}
\def\hom{\mathrm{Hom}}

\def\rela{{\rm Z} \kern -0.17cm \raise 0.05cm\hbox{\tiny{\it /}}\kern
  0.05cm}
\def\petitrela{{\rm Z} \kern -0.11cm \raise 0.05cm\hbox{{\mini /}}}
\def\ratio{{\rm Q} \kern -0.216cm\rule{0.25pt}{8pt}\kern 0.216cm}
\def\petitratio{{\rm Q} \kern -0.16cm\rule{0.15pt}{5pt}\kern 0.13cm}
\def\complex{{\rm C} \kern -0.20cm\rule{0.25pt}{7.5pt}\kern 0.20cm}
\def\petitcomplex{{\rm C} \kern -0.13cm\rule{0.15pt}{5pt}\kern 0.13cm}
\newcommand{\PP}{{\rm I} \! {\rm P}}

\newcommand{\OO}{{\cal O}}

\newlength{\elargir}
\title{}
\author{Han Fr\'ed\'eric}
\begin{document}
\phantom{.}
\vskip1cm
\noindent{\bf \LARGE Geometry of  genus 9 Fano 4-fold}
\vskip1cm
\noindent{\bf Han Fr\'ed\'eric}

{\small Institut de Math\'ematiques de Jussieu,
  
  Universit\'e Paris 7 - Case Postale 7012,
  
  2 Place Jussieu
  
  75251 Paris Cedex 05  France

  email: han@math.jussieu.fr}

\bigskip
\noindent{\large \bf Introduction:}

\def\Gw{G_\omega }
\def\Pw{P_\omega }
\def\PW{\PP(W)}
\def\EE{{\cal E}}
\def\II{{\cal I}}

Let $W$ be a $6$ dimensional vector space over the complex number,
equiped with a non degenerate symplectic form $\omega$. Let $\Gw$ be
the grassmannian of $\omega$-isotropic $3$-dimensional vector subspaces
of $W$. Considering the Plucker embedding,  the intersection of $\Gw$ with a generic
codimension $2$ linear subspace is the Mukai model of a smooth Fano manifold of dimension
$4$, 
genus $9$ index $2$ and picard number $1$. We will also note by $\Pw$
the $13$ dimensional projective space spanned by $\Gw$ under it's
plucker embedding.
\begin{description}
\item[\sc Notations:]  \ \\
  {\it In all the paper, $B$ will be a general double hyperplane section of $\Gw$. For a hyperplane $H$ of $\Pw$, we define  $\bar{H}=H\cap
\Gw$, and for any $u\in \Gw$, the corresponding plane of $\PP(W)$ will be noted
$\pi_u$.} 
\end{description}
\begin{abstract}
On a genus $9$ Fano variety, Mukai's construction gives a natural rank $3$ vector
bundle, but curiously in dimension $4$, another phenomena appears. In the 
first part of this article, we will explain how to construct on a
Fano $4$-fold of genus $9$ (named $B$), a canonical set of four stable  
vector bundles of rank $2$, and prove that they are rigid. Those
bundles was already known to A. Iliev and K. Ranestad  in \cite{IR}
and the results of this section could be viewed as particular cases of
the work of A. Kuznetsov (Cf \cite{K}). We'd like here to show their
consequences in the geometry of the $4$-fold, and study Zak duality in
this case.

Indeed, this ``four-ality'' (Cf \cite{M}) is also present in the geometry of
lines included in $B$, and also in 
the Chow ring of $B$. In section 2  we show that the variety of lines in $B$, is 
an hyperplane section of $\PP_1\times\PP_1\times\PP_1\times\PP_1$. This description is
explicit and  could also be interesting in terms of Freudenthal
geometries. Then in section 3, we
compute the  Chow ring of $B$ which appears to have a rich structure in codimension $2$.

The $4$ bundles constructed can 
embed $B$ in a Grassmannian $G(2,6)$, and the link with the order one congruence
discovered by E. Mezzeti and P de Poi in \cite{MdP} will be done in section 4. In
particular we will prove that the generic fano variety of genus $9$ and dimension $4$ can
be obtained by their construction, and explain the choices involved. We will also describe 
in this part the normalization of the non quadraticaly normal variety  they constructed,
and also its variety of plane cubics.
\end{abstract}
\begin{description}
\item[\bf Aknowledgements:] \ \\
  I'd like to thanks L. Gruson for his constant interest on this work, and also
  F. Zak and C. Peskine for fruitful discussions. I'd like also to
  thanks K. Ranestad and A. Kuznetsov for pointing and explaining their works.
\end{description}
\section{Construction of rank 2 vector bundles on $B$}\label{fibres}
This part is devoted to the construction of a canonical set of $4$ stable and rigid rank
$2$ vector bundles on $B$. Many results are done in a universal way in
\cite{K}, but we detail this description to use it in next sections. Let's first
recall some classical geometric properties of $\Gw$. ( Cf \cite{I})
The union of the tangent spaces to $\Gw$ is a quartic hypersurface of
$\Pw$, so a general line of $\Pw$ has naturally $4$ marked
points. Dually, as the variety  $B$ is given by a pencil $L$ of hyperplane sections of
$\Gw$, there are in this pencil, $4$ hyperplanes $H_1,\dots,H_4$
  tangent to $\Gw$. Denoting by $u_i$ the contact point of $H_i$ with $\Gw$, we
  will first construct a rank $2$  sheaf on $H_i\cap \Gw$ with singular locus
  $u_i$, and it's restriction to $B$ will be the vector bundle. 

\subsection{Data  associated to a tangent hyperplane section}\label{sectionhyp}
Let $u\in\Gw$, and  $H$ be a general hyperplane tangent to $\Gw$ at $u$.  For
any $v$ in $\Gw$, denote by $\pi_v$ the corresponding projective subspace of
$\PW$, and consider the hyperplane section of $\Gw$:
$$\bar{H}_u=\{v\in \Gw, \pi_v\cap \pi_u \neq \emptyset\}$$
It's proved in \cite{I} the following:
\begin{lemma}\label{LC}
There is a conic $C$ in $\pi_u$ such that $v\in H\cap \bar{H}_u \iff \pi_v\cap C
\neq \emptyset$. For $H$ general containing the tangent space of $\Gw$ at $u$,
$C$ is smooth. Furthermore, $H\cap \bar{H}_u$ contains the tangent cone
$T_u\Gw\cap \Gw=\{v \in \Gw | \dim \pi_v \cap 
\pi_u \geq 1\}$ which is embedded in $\Pw$ as a cone over a veronese surface.
\end{lemma}
Let $Z$ be the following incidence:
$$Z_H=\{(p,v)\in C\times (H\cap \Gw) | p \in \pi_v \}$$ 
Identifying $C$ with $\PP_1$, we denote by $q_1$ and $q_2$ the projections from
 $C\times \Gw$ to $\PP_1$ and to $\Gw$, and  by $L$ the $SL_2$-representation $H^0
\OO_{\PP_1}(1)$. Restricting the surjection  $L\otimes \OO_{\bar{H}} \to
q_{2*}\OO_{Z_H}(1,0)$ to the hyperplane section $\bar{H}$, we obtain:
\begin{proposition}\label{defE}
  The sheaf $\EE$ defined by the following exact sequence:
  $$0 \to \EE \to L\otimes \OO_{\bar{H}} \to q_{2*}\OO_{Z_H}(1,0) \to 0$$
  is reflexive of rank $2$, $c_1(\EE)=-1$ and is locally free outside $u$. 
\end{proposition}
\proof  From lemma \ref{LC}, the support of $q_{2*}\OO_{Z_H}(1,0)$ is $H\cap \bar{H_u}$
so it is an hyperplane section of $\bar{H}$, hence $\EE$ has rank $2$ and
$c_1(\EE)=-1$. Furthermore, for $v$ in $q_2(Z_H)$, the fiber of the restriction of $q_2$
to $Z_H$: $q_{2|Z_H}^{-1}(v)$ has
length $1$ if $v$ is not in $\bar{H}_u$, it has lenght $2$ if $v
\in\bar{H}_u-\{u\}$ and $q_{2|Z_H}^{-1}(u)$ is the curve $C$. As $Z_H$ and
$\bar{H}-u$ are   smooth, the
sheaf $q_{2*}\OO_{Z_H}(1,0)$ has projective dimension $1$ outside $u$, hence
$\EE$ is locally free outside $u$.\fin

\bigskip
Denote by $S_iL$ the
$SL_2$-representation $H^0(\OO_{\PP_1}(i))$ and by $K$ and $Q$ the tautological
bundles on $\Gw$, such that the following sequence is exact.
$$0\to K\rightarrow W\otimes \OO_{\Gw} \rightarrow Q\to 0$$
\begin{proposition}\label{resq2*}
  For $i>0$ we have $R^iq_{2*}\OO_{Z_H}(1,0)=0$, and the resolution of
  $q_{2*}\OO_{Z_H}(1,0)$ in $\Gw$ is given by the following   exact sequence:
  \begin{equation}
    \label{eqresq2*}
  0 \to S_3L \otimes\OO_{\Gw}(-1) \to L\otimes \bigwedge^2 Q\dual \to L\otimes \OO_{\Gw}
  \to q_{2*}\OO_{Z_H}(1,0) \to 0   
  \end{equation}
\end{proposition}
\proof
We consider the injection from $q_1^*(\OO_{\PP_1}(-2))$ to $W\otimes
\OO_{\PP_1\times \Gw}$ given by the conic $C$. So the incidence $Z_H$ is the
locus where the map from $q_1^*(\OO_{\PP_1}(-2)) \oplus q_2^*K$ to $W\otimes
\OO_{\PP_1\times \Gw}$ is not injective, hence  $Z_H$ is obtained in $\PP_1
\times \Gw$ as the zero locus of a section of the bundle
$\OO_{\PP_1}(2)\boxtimes Q$. 

Let ${\cal K}^{.}$ be the Koszul complex
$\bigwedge\limits^i(\OO_{\PP_1}(-2)\boxtimes Q\dual)$ of this
section. We obtain the proposition \ref{resq2*} from the Leray
spectral sequence applied to ${\cal K}^.$ twisted by $\OO_{\PP_1
  \times \Gw}(1,0)$. \fin

Furthermore, we deduce from Bott's theorem on $\Gw$ the following:
\begin{corollary}
   We have the following equality
   $L=H^0(\OO_{Z_H}(1,0)=H^0(q_{2*}\OO_{Z_H}(1,0))$,  and for $i>0$, all
   the groups $H^i(\OO_{Z_H}(1,0)$  and $H^i(q_{2*}\OO_{Z_H}(1,0))$
   are zero. For $i\geq 0$ all the groups $H^i(q_{2*}\OO_{Z_H}(1,-1))$
   and $H^i(q_{2*}\OO_{Z_H}(1,-1))$ are zero.
 \end{corollary}
 \proof
 We will prove that on the isotropic Grassmannian $\Gw$, the bundles
$\bigwedge\limits^i Q\dual$ and $(\bigwedge\limits^i Q\dual)(-1)$ are acyclic
for $i\in\{1,2,3\}$. Indeed, with the notations of \cite{W} 4.3.3 and 4.3.4, they
correspond to the partitions $(0,0,-1)$, $(0,-1,-1)$, $(-1,-1,-1)$, $(-1,-1,-2)$,
$(-1,-2,-2)$, 
$(-2,-2,-2)$. Now recall that the half sum of positive roots is $\rho=(3,2,1)$, so
$\alpha+\rho$ 
either contains a $0$ or is $(2,1,-1)$. So in all cases $\alpha+\rho$ is invariant by a
signed permutation, and the sheaves are acyclic. The corollary is now a direct consequence
of this acyclicity. \fin  

 \begin{corollary}\label{acyE}
The sheaves $\EE$ and $\EE(-1)$ are acyclic. The vector space
$V=H^0(\EE(1))$ has dimension 6 and $\forall i>0, h^i(\EE(1))=0$, and
$\EE(1)$ is generated by its global sections. 
\end{corollary}
\proof The acyclicity of $\EE$ and $\EE(-1)$ is a direct consequence
of  the definition of $\EE$ and of the previous corrolary.

To obtain the second assertion, we restrict the sequence
\ref{eqresq2*} to the hyperplane section $\bar{H}$, so we obtain the
following monad\footnote{A monad is complex exact at all terms
  different from the middle one. }:
$$0 \to S_3L \otimes\OO_{\bar{H}}(-1) \to L\otimes \bigwedge^2 Q\dual_{\bar{H}} \to
\EE \to 0  $$
whose cohomology is $\Tor^1(q_{2*}(\OO_{Z_H}(1,0)),
\OO_{\bar{H}})$ which is equal to $q_{2*}(\OO_{Z_H}(1,-1))$ because
$Z_H \subset q_2^{-1}(H)$. Twisting this monad by
$\OO_{\bar{H}}(1)$ we obtain that $H^0(\EE(1))$ is the quotient of
$L\otimes W$ by $S_3L\oplus L$ because
$W=H^0(Q_{\bar{H}})$. Furthermore, the right part of the monad gives a
sujection from $L\otimes Q_{\bar{H}}$ to $\EE(1)$. But $L\otimes Q_{\bar{H}}$ is
globally generated, so $\EE(1)$ is also generated by its global sections.

The vanishing of $h^i(\EE(1))$ for $i>0$ is a corollary of the vanishing
of $h^i(q_{2*}(\OO_{Z_H}(1,0))$, $h^i(Q_{\bar{H}})$ and
$h^i(\OO_{\bar{H}})$ for $i>0$. \fin

We can remark that the two vector spaces $V$ and $W$ of dimension 6
have not the same role. More precisely, the conic $C$ gives a marked
subspace of $W$ so that we have the following:
\begin{remark}\label{SL2WV}
  The  tangent hyperplane $H$ gives canonically the $SL_2$-equivariant
  sequences:
  \begin{center}
  $0 \to S_2L \to W \to S_2L\to 0$ and $0\to L \to V \to S_3L \to 0$    
  \end{center}
\end{remark}
\subsection{The 4 rank 2 vector bundles on $B$}
The pencil of hyperplanes defining $B$ contains the $4$ tangent hyperplanes
$H_i$, so we can apply the previous construction to construct a rank $2$ sheaf
$\EE_i$ on each of the $\bar{H}_i$, and define by $E_i$ the restriction of $\EE_i$ to
$B$. Because $B$ is smooth, it doesn't contain the contact points $u_i$, so
$E_i$ is locally free on $B$.
\begin{corollary}\label{acyEi}
All the cohomolgy groups of the vector bundles $E_i$  vanish. In
particular, the rank $2$ vector bundles $E_i$ are stable. The
vector space $H^0(E_i(1))$ has dimension $6$, and $\forall j>0,
h^j(E_i(1))=0$. The bundles $E_i(1)$ are generated by their global sections.
\end{corollary}
\proof It's a direct consequence of corollary \ref{acyE}, because  $B$ is a
hyperplane section of $\bar{H_i}$. (Note that the stablility condition
for a $E_i$ is equivalent to $h^0E_i=0$) \fin
\subsection{The restricted incidences}
Now, for each of the $4$ hyperplanes $H_i$ containing $B$ and tangent to
$\Gw$ at some point $u_i$, let $C_i$ be the conic of the projective plane $\pi_{u_i}$
constructed in \ref{sectionhyp}.  Consider the restriction
of the incidences $Z_{H_i}$ to $B$. In other words, let $Z_i,Z'_i$ be:
$$Z_i=\{(p,v)\in C_i \times B | p \in \pi_v \},\ \ Z'_i=\{(p,v)\in Z_i | \dim
(\pi_v\cap \pi_{u_i})>0 \}$$
where $q_1$ and $q_2$ still denote the projections from $C_i\times \Gw$ to $C_i$ and $\Gw$.
\begin{remark}\label{QetV}
 Let $p$ be a fixed point of $C_i$. The scheme $Z_{i,p}=q_2(q_1^{-1}(p)\cap
 Z_i)$ is a 2 dimensional irreducible quadric  in 
 $\Pw$. The restriction of $q_2$ to $Z_i'$ is a double cover of a veronese
 surface $V_i=q_2(Z_i')$.
\end{remark}
\proof In fact $\{b\in \Gw | p \in \pi_v \}$ is a smooth quadric of dimension 3 (Cf
\cite{I}), so it doesn't contain planes. This scheme is included in $H_i$, so
$Z_{i,p}$ is just an hyperplane section of this smooth quadric. It is also
proved in \cite{I} that $\{v\in \Gw | \dim \pi_{u_i}\cap \pi_v>0 \}$ is a c\^one
over a veronese surface of vertex $u_i$. As $u_i\notin
B$, the surface $V_i$ is the intersection of this cone with an  hyperplane which doesn't
contain the vertex $u_i$, so it's a veronese
surface. \fin
\begin{description}
\item[\sc Notations:] \ \\
{\it We denote by $\sigma_i$  the class of a point on
$C_i$, and $h_3$ the class of a hyperplane in $\Pw$. (the plucker
embedding of $\Gw$).}
\end{description}
\begin{proposition}\label{modelZi}
  The incidence $Z_i$ is a divisor of classe $2h$ in
  $$\Pi=Proj(\OO_{\PP_1}(2)\oplus S_2L\otimes \OO_{\PP_1})$$
  Furthermore  we have $h_3\sim h+2\sigma_i$,  where  $h\sim
  \OO_{\Pi}(1)$ and $\sigma_i$ is also the class of a point on the
  base of the fibration $\Pi$ . The divisor $Z'_i$ of $Z_i$ is equivalent to
  $h-2\sigma_i$. 
\end{proposition}
\proof  Denotes by $e_i$ the image of the map from $\OO_{\PP_1}(-2)$
to $W\otimes \OO_{\PP_1}$ associated to $C_i$. Choose an element $\phi'$ of
$\bigwedge\limits^3 W\dual$ such that  $\ker \phi'$ gives an hyperplane section
of $\Gw$ containing $B$ and  different from the $\bar{H}_i$. (i.e $\phi'(u_i)\neq
0)$. We can remark that the 
incidence $Z_i$ is 
given over $\PP_1$ by the isotropic $2$-dimensional subsbaces $l$ of
$e_i^{\bot}\over e_i$, such that $\phi'(e_i\otimes \wedge^2 l)=0$,
because the condition $\phi_i(e_i\otimes \wedge^2 l))=0$ is already
satisfied by the definition of $C_i$ and lemma \ref{LC}. (where $\phi_i$ denotes
a trilinear form of kernel $H_i$) 

The bundle $e_i^{\bot}$ is isomorphic to $S_2L\otimes\OO_{\PP_1} \oplus L\otimes
\OO_{\PP_1}(-1)$ where the trivial factor $S_2L$ correspond to the
plane $\pi_{u_i}$. So\footnote{As $SL_2$-represention, we will
  identify $L$ with its dual.} the bundle $e_i^{\bot}\over e_i$ is 
isomorphic to $L(1) \oplus L(-1)$ where those factors are isotropic for the
symplectic form induced by $\omega$. We can take local basis $s_0,s_1$ and
$s_2,s_3$ of each factors such that the form induced by $\omega$ is $p_{0,2}+p_{1,3}$
where  $p_{i,j}$ denotes the Plucker coordinates associated to the $s_i$.

So the relative isotropic grassmannian $\Gw(2,{e_i^{\bot}\over e_i})$ is the
intersection of $G(2,{e_i^{\bot}\over e_i})$ with the subsheaf
$\OO_{\PP_1}(2)\oplus \OO_{\PP_1}(-2)+S_2L\otimes\OO_{\PP_1}$ where the factor
$\OO_{\PP_1}(2)$ still correspond to $s_0\wedge s_1$. 

Now we need to compute the kernel of the map  $e_i\otimes
\bigwedge\limits^2({e_i^{\bot}\over e_i})\stackrel{\phi'}{\rightarrow}
\OO_{\PP_1}$. But the assumption $\phi'(u_i)\neq 0$ proves that it is
$\OO_{\PP_1}(-4)\oplus L\otimes L\otimes\OO_{\PP_1} $.

So we have an exact sequence:
$$0 \to \OO_{\PP_1}(-2)\oplus S_2L\otimes \OO_{\PP_1}\to
\bigwedge\limits^2({e_i^{\bot}\over e_i}) \stackrel{( ^\omega
  _{\phi'})}{\longrightarrow} \OO_{\PP_1}\oplus e_i\dual \to 0$$
and $Z_i$ is a divisor of class $2h$ in $Proj(\OO_{\PP_1}(2)\oplus S_2L\otimes
\OO_{\PP_1})$. The 
relation $h_3\sim h+2\sigma_i$ is given by the map $e_i\otimes
\bigwedge\limits^2({e_i^{\bot}\over e_i}) \to \bigwedge\limits^3 W$.

The divisor $Z'_i$ of $Z_i$ is locally given by the vanishing of the exterior product
with $s_0\wedge s_1$ so it equivalent to $h-2\sigma_i$. \fin 

\bigskip
We will now study the relation between the conormal bundle of $Z_i$ in
$\PP_1\times B$ and the bundle $E_i$.
\subsection{Deformations of $E_i$}
\begin{lemma}\label{Brel}
  We have the following exact sequence:
  $$0 \to \OO_{\PP_1\times B}(-\sigma_i -h_3) \to q_2^* E_i  \to \II_{Z_i}(\sigma_i)\to
  q_2^*(R^1q_{2*}\II_{Z_i})(-\sigma_i) \to 0   $$
  (where $q_1$ and $q_2$ denotes the projections from $\PP_1\times B$ to $\PP_1$ and $B$)
\end{lemma}
\proof From the resolution of the diagonal of $\PP_1\times \PP_1$, we obtain the
relative Beilinson's spectral sequence:
$$E^1_{a,b}=(\bigwedge\limits^{-a}\omega_{q_2}(\sigma_i))\otimes
R^{b}q_{2*}(\II_{Z_i}((1+a)\sigma_i)) \Longrightarrow \II_{Z_i}(\sigma_i)$$
By the definition of $E_i$ (Cf prop
\ref{defE}) we have $E_i=q_{2*}(\II_{Z_i}(\sigma_i))$. Furthermore, the projection $q_2(Z_i)$
is an hyperplane section of $B$, so  $q_{2*}\II_{Z_i}=\OO_B(-1)$. We can
conclude, remarking that $R^1q_{2*}(\II_{Z_i}(\sigma_i))=0 $, because the
restriction of $q_2$: $Z_i \stackrel{q_{2|Z_i}}{\longrightarrow} q_2(Z_i)$ has all its fibers of length at most $2$. \fin

NB: The support of $R^1q_{2*}\II_{Z_i}$ is the natural scheme structure (Cf
\cite{GP}) on the scheme of fibers of $q_2$ intersecting $Z_i$ in length $2$ or
more. It is the veronese surface $V_i=q_2(Z'_i)$.

So the previous lemma can now be translated in the following:
\begin{corollary}\label{zeroL}
  The scheme $q_2^{-1}(V_i)\cup Z_i$ is in $\PP_1\times B$ the zero locus of a section of
  the bundle $q_2^*   E_i(1,1)$. 
\end{corollary}
This gives also a geometric description of the marked pencil of sections of
$E_i(h_3)$ given by the natural inclusion $L\subset V$ found in remak
\ref{SL2WV}. Indeed, if we fixe a point $p$ on $C_i$, the restriction to
$q_1^{-1}(p)$ of the section obtained in corollary \ref{zeroL} gives with the
notations of lemma 
\ref{QetV} the following: 
\begin{corollary}\label{pinceau}
  For any point $p$ on the conic $C_i$, the vector bundle $E_i(h_3)$ has a
  section vanishing on $Z_{i,p}\cup V_i$.  
\end{corollary}
We can now study the restriction of $E_i$ to $Z_i$.
\begin{proposition}\label{EiZi}
  The restriction $E_{i|Z_i}$ of the vector bundle $q_2^*E_i$ to $Z_i$ fits into the
  following exact sequence: 
  $$0 \to \OO_{Z_i}(h_3-3\sigma_i) \to q_2^*E_{i|Z_i}(h_3) \to \OO_{Z_i}(3\sigma_i)
  \to 0 $$  
\end{proposition}
\proof Fix a point $p$ on $C_i$, and consider the corresponding section of
$E_i(h_3)$ constructed in corollary \ref{pinceau}. Its pull back gives a section
of $q_2^*E_i(h_3)$ vanishing on $q_2^{-1}(Z_{i,p}\cup V_i)$, so its
restriction to $Z_i$ gives a
section of $q_2^*E_{i|Z_i}(h_3-\sigma_i-Z_i')$. Now,
using the  computation  of the class of $Z_i'$ in $Z_i$ made in  proposition
\ref{modelZi}, namely  that $\OO_{Z_i}(Z_i')$ is  $\OO_{Z_i}(h_3-4\sigma_i)$, it
gives a section of $E_{i|Z_i}(3\sigma_i)$. We have to prove that it is a non
vanishing section. To obtain this, we compute the second Chern's class of
$E_{i|Z_i}(3\sigma_i)$. We will show that its image in the Chow ring of $\PP_1
\times B$ is zero. Denote by $a_i$ the second Chern class of $E_i$. From
the lemma \ref{zeroL}, we obtain the class of $Z_i$ in $\PP_1 \times B$:
$[Z_i]=a_i+h_3.\sigma_i -[V_i]$. So we can compute
$[Z_i].c_2(E_i(3\sigma_i))$. It is $(a_i+h_3.\sigma_i
-[V_i]).(a_i-3h_3\sigma_i)$, but we will compute in proposition \ref{ChowB} the
Chow ring of $B$, and this class vanish. \fin
%
\begin{corollary}
  The vector bundles $E_i$ are  rigid, in other words we
  have $\ext^1(E_i,E_i)=0$.  
\end{corollary}
\proof From the corollary \ref{zeroL} we have an exact sequence on
$\PP_1\times B$:
$$ 0 \to q_2^* E_i(-\sigma_i) \to q_2^*(E_i)\otimes q_2^*(E_i)(h_3)
\to q_2^*E_i(h_3+\sigma_i) \to (q_2^*E_i(h_3+\sigma_i))_{|Z_i \cup
  q_2^{-1}(V_i)} \to 0$$
The bundle $q_2^* E_i(-\sigma_i)$ is acyclic, and  the corollary
\ref{acyEi} gives $H^0(q_2^*E_i(h_3+\sigma_i))=L\otimes V$ and
$H^1(q_2^*E_i(h_3+\sigma_i))=0$. The liaison exact sequence twisted by
$q_2^*(E_i(h_3+\sigma))$ is:
$$ 0 \to q_2^*E_i(h_3)\otimes \OO_{q_2^{-1}(V_i)}(\sigma_i-Z_i') \to
q_2^*E_i(h_3+\sigma_i))_{|Z_i \cup  q_2^{-1}(V_i)} \to E_{i|Z_i}(h_3+\sigma_i)\to 0$$
As $\sigma_i -Z_i'$ have degree $-1$ along the fibers of $q_2:
q_2^{-1}(V_i) \to V_i$, the bundle $q_2^*E_i(h_3)\otimes
\OO_{q_2^{-1}(V_i)}(\sigma_i-Z_i')$ is acyclic, so the cohomology of
$q_2^*E_i(h_3+\sigma_i))_{|Z_i \cup  q_2^{-1}(V_i)} $ can be computed
with its restriction to $Z_i$. The propositions \ref{EiZi} and
\ref{modelZi} show that $H^0 E_{i|Z_i}(h_3+\sigma_i)=S_2L\oplus
S_2L\oplus S_4L$. In conclusion, we have the exact sequence:
$$0 \to \hom(E_i,E_i) \to L\otimes V \to S_2L\oplus S_2L \oplus S_4L
\to \ext^1(E_i,E_i)  \to 0$$
By the corollary \ref{acyEi}, the bundle $E_i$ is stable, so it is
simple, in other words we have $\hom(E_i,E_i)=\complex $, and the
above exact sequence gives $\ext^1(E_i,E_i)=0$. \fin

\section{The variety of lines in $B$}

\begin{remark}\label{dtesGw}
  Let $\delta$ be an isotropic line of $\PW$. The set of isotropic planes of
  $\delta^{\bot}$ containg $\delta$ form a line in $\Gw(3,6)$, and all the line
  in $\Gw(3,6)$ are of this type for a unique element of
  $\Gw(2,6)$. In other words, the   variety of lines in $\Gw(3,6)$ is
  naturally $\Gw(2,6)$.  
\end{remark}
\def\KB{{(\frac{K_2^{\bot}}{K_2})}\dual(h_2)}%
\begin{description}\label{notations} 
\item[\sc Notations:] \ \\
  {\it A point  of $\Gw(2,6)$ will be denoted by a
minuscule letter, and the corresponding line in $\Gw(3,6)$ by the
majuscule letter. The  variety of lines included in $B$ will be noted
$F_B$. Denote by $I$ the incidence point/line in
$B$. In other words:
$$ I=Proj(\KB)=\{(\delta,p)\in F_B\times B | p\in \Delta \}\subset \Gw(2,6)\times
\Gw(3,6)$$
The projections from $I$ to $F_B$ and $B$ will be denoted by $p_1$ and $p_2$.}
\end{description}
\subsection{A morphism from $F_B$ to   $\PP_1\times\PP_1\times\PP_1\times\PP_1 $}
Each of the $4$ conics $C_i$ will enable us to construct a morphism
from $F_B$ to $\PP_1$. We have the following  geometric hint to expect
 at least a rational map: A general element $\delta$ of $F_B$ gives an isotropic
$2$-dimensional subspace $L_\delta$ of $W$.  In general, the projectivisation of
$L_\delta^{\bot}$  meets the plane containg $C_i$ in a point
$p$. There is at least an element $m$ of $\Delta$ such that
$p\in\pi_m$, so $m\in\Delta\subset B\subset H_i$. Now, the definition of
$H_i$ and lemma \ref{LC} prove that $p$ must be on $C_i$. 

But to show that it is everywhere defined, we will use the vector
bundle $E_i$. We start by constructing line bundles on $F_B$.
\begin{remark}\label{secVi}
  Any line $\Delta$ included in the hyperplane section $\bar{H}_{u_i}=q_2(Z_i)$ of $B$
  intersect the veronese surface $V_i$. Furthermore, any such line is in a  quadric
  $Z_{i,p_i}$ for a unique point   $p_i$ of $C_i$. So the set $v_i=\{\delta \in F_B |
  \Delta \subset \bar{H}_{u_i}\}$ is a divisor in $F_B$.
\end{remark}
\proof As the line $\Delta$ is in $\bar{H}_{u_i}$, we have from lemma \ref{LC} that for
any $b\in \Delta$, the plane $\pi_b$ intersect the conic $C_i$. The line
$\PP(L_{\delta})$ must intersect $C_i$ in some point $p_i$. Indeed, if it was not the case, this line would be
orthogonal to $C_i$, so it would be in the plane $\pi_{u_i}$, but any line in this plane
intersect $C_i$. 

So the line $\Delta$ is in the quadric $Z_{i,p_i}$. Note that $\PP(L_{\delta})\cap C_i$
can't contain another point,  because it would imply that $\Delta \subset
V_i$. Furthermore, the proposition \ref{modelZi} implies that $Z_{i,p_i}\cap V_i$ is a
plane section of the quadric $Z_{i,p_i}$, so $\Delta$ intersect $V_i$ in a single point. \fin
\begin{remark}
For any point $p_i$ of $C_i$, the scheme $p_2^{-1}(Z_{i,p_i})$ has several irreducible
components of dimension $2$,  but  some of these component are contracted by $p_1$ to a
curve. Denote by $A_{i,p_i}$ the $2$ dimensional part of  $p_1(p_2^{-1}(Z_{i,p_i}))$. 
\end{remark}
\proof The components of $p_2^{-1}(Z_{i,p_i})$ corresponding to the lines included in
$Z_{i,p_i}$ are contracted to  curves. \fin
\begin{proposition}\label{alphai}
The sheaf  $p_{1*}p_2^* E_i$ is a line bundle on $F_B$. There is a natural map $f_i$ from 
$H^0(\OO_{C_i}(\sigma_i))\dual\otimes \OO_{F_B}$ to the dual bundle of
$(p_{1*}p_2^* E_i)$. The image of $f_i$ is also a line bundle on $F_B$, we will  denote it  by
$\OO_{F_B}(\alpha_i)$. By construction, for any $p_i\in C_i$, the divisor 
$A_{i,p_i}$ will be in the linear system $|\OO_{F_B}(\alpha_i) |$.
\end{proposition}
\proof By the corollary \ref{acyEi}, the bundle $E_i$ is a quotient of
$6\OO_B(-1)$ and by  definition \ref{defE}, it is a subsheaf of
$2\OO_B$. So its restriction to any line $\Delta$ included in $B$ must
be $\OO_{\Delta}\oplus \OO_{\Delta}(-1)$. So $R^1p_{1*}p_2^*E_i=0$ and
$p_{1*}p_2^*E_i$ is a line bundle.
Denote this line bundle by $\OO_{F_B}(-\alpha_i')$.
Dualising and twisting the exact sequence defining $E_i$, we have the following exact
sequence:
\begin{equation}
  \label{eq:Li}
0 \to  L\otimes \OO_B(-2h_3) \to E_i(-h_3) \to \LL_i \to 0
\end{equation}
where $\LL_i$ is supported on the hyperplane section $\bar{H}_{u_i}$, and is singular along
the veronese surface $V_i$. As the incidence $I$ is $Proj(\KB)$ (where $K_2$ is the
tautological subbundle of $W\otimes \OO_{\Gw(2,W)}$), the relative dualising
sheaf $\omega_{p_1}$ is $\OO_I(2h_2-2h_3)$. So we have
$R^1p_{1*}(p_2^*E_i(-h_3))=\OO_{F_B}(\alpha_i'-2h_2)$. So the base locus of this pencil of
sections of $\OO_{F_B}(\alpha_i')$ is the support of $R^1 p_{1*}p^*_2(\LL_i)$.  We will
now  prove that this sheaf is a line bundle on $v_i$. 

The morphism $p_1$ is projective of relative dimension
$1$. So,  for any point $\delta$ of $F_B$ and any coherent sheaf $F$ on $B$,  the fiber
$(R^1p_{1*}p_2^* F)_{\delta}$ is 
$H^1(F\otimes \OO_\Delta)$, where $\Delta$ is the line in $B$ corresponding to
$\delta$. The restriction of the sequence \ref{eq:Li} to $\Delta$ gives the surjection:
\begin{equation}
  \label{restricDelta}
 2\OO_\Delta(-2) \to \OO_\Delta(-1) \oplus \OO_\Delta(-2) \to \LL_i\otimes \OO_\Delta
\to 0  
\end{equation}
When the line $\Delta$ is not in the hyperplane section $\bar{H}_{u_i}$, the sheaf
$\LL_i\otimes \OO_\Delta$ is supported by the point $\bar{H}_{u_i}\cap
\Delta$, so  in this case we have $h^1(\LL_i\otimes \OO_\Delta)=0$. Now, when the line
$\Delta$ is in  $\bar{H}_{u_i}$, the sheaf $\LL_i\otimes \OO_\Delta$ has generic rank $1$
because the veronese surface $V_i$ can't contain the line $\Delta$. We have proved in
lemma \ref{secVi} that $\Delta$ intersect $V_i$, hence for any element of $L$, the
corresponding section of $E_i(h_3)$ vanishes on $\Delta$, so the map $2\OO_\Delta(-2) \to
\OO_\Delta(-2)$ induced by the  sequence (\ref{restricDelta}) is zero, and for any $\delta$ in $v_i$, we
have $h^1(\LL_i\otimes \OO_\Delta)=1$, and  $R^1 p_{1*}p^*_2(\LL_i)$ is a line bundle on
$v_i$. 

So we have proved that the base locus of $|(p_{1*}(p^*_2 E_i))\dual|$ is the divisor
$v_i$. In other words, the image of $f_i$ is the line bundle $\OO_{F_B}(\alpha_i)=
(p_{1*}(p^*_2 E_i))\dual(-v_i) $, which is by construction a quotient of $L\otimes
\OO_{F_B}$. By definition $A_{i,p_i}$ is the closure of $\{\delta \in F_B | \mbox{lenght}(\Delta\cap
Z_{i,p_i})=1\}$ which was identified set theoretically with an element of the linear
system $|\alpha_i|$, so we conclude the proof with  lemma \label{reduit}:  \fin
\begin{lemma}\label{reduit}
  For a generic choice of a point $p_i$ on $C_i$,  the support of the sheaf 
 \linebreak$R^1p_{1*}(p_2^* I_{Z_{i,p_i}\cup V_i}(-h_3))$ represent the class  $\alpha'_i$, and  all
  its irreducible components are reduced.
\end{lemma}
\proof 
First notice that the point $p_i$ on $C_i$,  gives a section of $E_i(h_3)$, so an exact sequence:
$$0 \to \OO_B(-2h_2) \to E_i(-h_3) \to I_{Z_{i,p_i}\cup V_i}(-h_3)\to 0$$
which gives a section of $\OO_{F_B}(\alpha_i')$ vanishing on the support of the sheaf
\linebreak$R^1p_{1*}(p_2^*I_{Z_{i,p_i}\cup V_i}(-h_3))$. But this is the definition in \cite{GP} of the scheme
structure on the set of lines included in $B$ and intersecting $Z_{i,p_i}\cup V_i$. So to
show that this scheme structure is reduced on each component, we have to prove that
$Z_{i,p_i}$ and $V_i$ are not in the ramification of the morphism: $p_2:I \to B$. But a
general point $m$ on $Z_{i,p_i}\cap V_i$ is the intersection of $Z_{i,p_i}$ and another
quadric $Z_{i,p'_i}$, so there pass  $4$ distinct lines through $m$, and from remark
\ref{secVi} there are no other lines in $B$ through $m$. So $m$ is not in the
ramification of the morphism $p_2:I\to B$ because it is of degree $4$. (Cf lemma
\ref{quasifini}). \fin 
\subsection{Description of the morphism}
In the previous section, we have contructed $4$ morphism $f_i$ from $F_B$ to
$\PP_1$. We will now prove that the morphism from $F_B$ to $\PP_1\times
\PP_1\times \PP_1\times \PP_1$ is an embedding for a generic $B$, and that its
image is an  hyperplane section.

\bigskip
\begin{description}
\item[\sc Notations:] \ \\
  {\it For a point $p_i$ in the conic $C_i$, we denote by ${\cal C}_B$
the affine cone over $B$. We consider:} 
$$F_{p_i}=\{\delta \in \Gw(2,W) | \delta \wedge p_i \in {\cal C}_B \}$$
\end{description}
Unfortunately, we have to remark that $f_i^{-1}(p_i)$ is not exactly $F_{p_i}\cap F_B$:
\begin{remark}\label{warningFpi}
  The intersection $F_{p_i}\cap F_B$ is equal to $\Gw(2,p_i^\bot)\cap F_B$. It contains
  $f_i^{-1}(p_i)$ and  residual curves corresponding to the lines included in
  the quadric $Z_{i,p_i}$.
\end{remark}
\proof If $\delta $ is already an element of $F_B$, then any isotropic plane of
$\PW$ containing the line $\PP(L_\delta)$ is an element of $B$ (Cf remark
\ref{dtesGw}). For $\delta \in \Gw(2,p_i^\bot)$, the plane 
 $\pi_{\delta \wedge p_i}$ is isotropic, so we have $\delta \wedge p_i \in {\cal C}_B$ and
 $\Gw(2,p_i^\bot)\cap F_B= F_{p_i}\cap
 F_B$. Now, remark that in  $\Gw(2,W)$ the scheme $\Gw(2,p_i^\bot)$ is the zero 
locus of a section of $K_2\dual$. The restriction of this section to $F_B$
vanishes on the 
divisor $A_{i,p_i}$ defined in proposition \ref{alphai}, so $F_B\cap F_{p_i}$ contains a divisor of class
$\alpha_i$  and the zero locus of a section of $(K_2\dual)_{|F_B}(-\alpha_i)$
corresponding to the lines includes in $Z_{i,p_i}$. \fin

\bigskip
Nevertheless, we have the following:
\begin{lemma}\label{Fij}
  For the generic double hyperplane section $B$ of $\Gw(3,W)$, we can find
  points $p_{i}$ (resp $p_{j}$) in the conics $C_i$ (resp $C_j$)
  such that $F_{p_{i}}\cap F_{p_{j}}$ is a smooth conic in $\Gw(2,W)$. Furthermore, this
  conic is in $F_B$ and represent the class $\alpha_i.\alpha_j$. 
\end{lemma}
\proof We can choose $p_i$ and $p_j$ respectively in $C_i$ and $C_j$ such that,
$p_i$ is not in $p_j^{\bot}$. Remark first that this implies that the intersection
$F_{p_i}\cap F_{p_j}$   is automatically included in $F_B$ because $p_i$ and $p_j$ are
never in the same isotropic plane. This also implies that the isotropic grassmannian 
$\Gw(2,p_i^{\bot}\cap p_j^{\bot})$ is a smooth $3$-dimensional quadric. The
intersection $F_{p_i}\cap F_{p_j}$ is given in this grassmannian by the $2$
conditions: $d\wedge p_i \in H_j$ and $d\wedge p_j \in H_i$ for an element $d$
of $\Gw(2,p_i^{\bot}\cap p_j^{\bot})$. Indeed, as $d\wedge p_i$ and $d\wedge
p_j$ represent isotropic planes, the conditions $d\wedge p_i \in H_i$ and
$d\wedge p_j \in H_j$ are automatically satisfied because $p_i\in C_i$ and
$p_j\in C_j$ (Cf lemma \ref{LC}). Let $\delta$ be the intersection of
$p_j^{\bot}$ and the $3$ dimensional vector space $U_i$ represented by the
contact point $u_i$ of $H_i$. The space $\delta$ is an element of
$\Gw(2,p_i^{\bot}\cap p_j^{\bot})$ such that $p_i\wedge \delta$ is not in
$H_j$ (because it is $u_i$,  and $B$ is smooth) but $p_j\wedge \delta \in H_i$
according to lemma \ref{LC}. So the two hyperplane sections are independent and
we have proved that $F_{p_i}\cap F_{p_j}$ is a (may be singular) conic.

Note that from the genericity of $B$ we could assume also that $\PP(p_i^{\bot}) \cap
C_j=\{a_1,a_2\}$ and $\PP(p_j^{\bot}) \cap C_i=\{b_1,b_2\}$ are $4$ distinct
points. According to lemma \ref{LC}, the lines $(p_i,a_1)$, $(p_i,a_2)$,
$(p_j,b_1)$, $(p_j,b_2)$ are in  $F_{p_i}\cap F_{p_j}$, but no three of those
lines  are in the same plane, so the conic $F_{p_i}\cap F_{p_j}$ contains $4$
points with no trisecant line. Hence the conic must be smooth. To conclude that this conic
represent the class $\alpha_i.\alpha_j$, we just have to prove that the residual curve of
$F_{p_i}\cap F_B$ don't intersect $F_{p_j}$. But if $\delta$ is such that $\Delta \subset
Z_{i,p_i}$, the line $\PP(L_\delta)$ contains the point $p_i$ which is not orthogonal to
$p_j$, so $\delta \notin F_{p_j}$. \fin
\begin{lemma}\label{Cijk}
  When $i,j,k$ are distinct, the morphism from $F_B$ to $C_i\times C_j \times
  C_k$ is  dominant.
\end{lemma}
\proof According to lemma \ref{Fij}, for a generic choice of $p_i\in C_i$ and
$p_j\in C_j$ the subvariety $F_{p_i}\cap F_{p_j}$ of $F_B$ is a smooth
conic, and we can also assume that the line $(p_ip_j)$ doesn't intersect
$C_k$. Assume that the induced map from  $F_{p_i}\cap F_{p_j}$ to $C_k$ is not 
dominant, then there is a point $p_k\in C_k$ such that 
$F_{p_i}\cap F_{p_j}\subset F_{p_k}$. So for any element $d$ of $F_{p_i}\cap
F_{p_j}$ the corresponding line $D$ is in the plane $\PP(<p_i,p_j,p_k>^{\bot})$. As
$p_k\notin (p_ip_j)$, the vector space  $<p_i,p_j,p_k>^{\bot}$ has dimension
$3$, and this contradicts the fact that $F_{p_i}\cap F_{p_j}$ is a smooth
conic. \fin

At this point, we need more details about the embedding of $F_B$ in $\Gw(2,W)$:

\begin{description}
\item[\sc Notations:] \ \\
  {\it Still denote by $K_2$ the tautological rank 2 subsheaf of $W\otimes
\OO_{\Gw(2,W)}$, and by $-h_2$ and $c_2$ its first and second Chern classes. As $\Gw(2,W)$
is an hyperplane section of $G(2,W)$, we will do the computations in $G(2,W)$. }
\end{description}
\begin{remark}\label{chowGw26}
 The Chow ring of $G(2,6)$ is
 $$\complex[h_2,c_2]/(h_2^5+3h_2c_2^2-4h_2^3c_2,-h_2^4c_2+3h_2^2c_2^2-c_2^3)$$
\end{remark}
We have:
\begin{lemma}\label{classFB}
 The variety $F_B$  is obtained in $ \Gw(2,W)$ as the zero locus of a
section of the bundle $\KB \oplus \KB$. For a general choice of $B$, $F_B$ is smooth with $\omega_{F_B}=\OO_{F_B}(-h_2)$. The class
of $F_B$ in the chow ring of $\Gw(2,W)$ is $4.(h_2^2-c_2)^2$, and $F_B$ has
degree $24$ in $\Gw(2,W)$. In the Chow ring of $F_B$, we have the extra relation: $h_2^2=3c_2$.
\end{lemma}
\proof To compute Chern classes, we just have to remark that $K_2^{\bot}$ is
isomorphic to the dual of the tautological 
quotient. The vanishing locus claim is a consequence of the fact that $B$ is a double
hyperplane section 
of $\Gw(3,W)$, and that we have from the definition of $I$ the equality: 
$p_{1*}(p_2^*\OO_{\Gw(3,W)}(h_3))=\KB$. (with the notations introduced in
{\S}\ref{notations}). So  the choice of a generic $2$-dimensional subspace of $H^0(
\OO_{\Gw(3,W)}(h_3)))$ correspond to the choice of a generic section of $\KB\oplus \KB$. 
Hence,  for a general choice of $B$, $F_B$ will be smooth because $K_2^{\bot}(h_2)$ is 
globally generated,  and so is $\frac{K_2^{\bot}}{K_2}(h_2)\simeq \KB$. We can
now conlude with the computations using the remark \ref{chowGw26}. \fin

\begin{lemma}\label{calculimage}
  In the Chow ring of $F_B$ we have $\alpha_i^2=0$, $\alpha_i.\alpha_j.h_2=2$ and
  $\alpha_i.\alpha_j.\alpha_k=1$ when $i,j,k$  are distinct. (where $h_2$ is the
  hyperplane section of $F_B$)
\end{lemma}
\proof In proposition \ref{alphai}  we already proved that  $\alpha_i^2=0$.
We have proved in the lemma \ref{Cijk}, that for a generic choice of
$p_i,p_j,p_k$ the intersection $F_{p_i}\cap F_{p_j}\cap F_{p_k}$ is not
empty. Furthermore, this intersection is included in the smooth conic
$F_{p_i}\cap F_{p_j}$ and the line $\Gw(2,p_i^{\bot}\cap p_j^{\bot} \cap
p_k^{\bot})$. So it must be a point because the intersection of a smooth conic
and a line in the $3$-dimensional smooth quadric $\Gw(2, p_i^{\bot}\cap
p_j^{\bot})$ must be empty or a point. So $\alpha_i.\alpha_j.\alpha_k=1$. \fin

\begin{lemma}\label{equiv}
The hyperplane section $h_2$ of $F_B$  is linearly
equivalent to $\alpha_1+   \alpha_2+   \alpha_3+   \alpha_4   $.
\end{lemma}
\proof 
We deduce from the lemma \ref{calculimage} that the image of the morphism from
$F_B$ to $\PP_1\times\PP_1\times\PP_1\times\PP_1$ is a divisor $\bar{F}$ of class
$(1,1,1,1)$. Furthermore, the lemmas \ref{Cijk} and \ref{calculimage} imply that $\psi: F_B \to
\bar{F}$ is a birationnal morphism.

Now we should notice that $\bar{F}$ is normal because  the restriction of  $\psi: F_B \to
\bar{F}$ to a double hyperplane section of $\PP_1\times\PP_1\times\PP_1\times\PP_1$ is not
ramified. Indeed, from the lemma \ref{calculimage} we have:
$(\alpha_1+\alpha_2+\alpha_3+\alpha_4)^3=24=h_2.(\alpha_1+\alpha_2+\alpha_3+\alpha_4)^2$,
so the restriction of $\psi$ to a general double hyperplane section of $\bar{F}$
is not ramified.

From the normality of  $\bar{F}$, we obtain a map $\psi^*(\omega_{\bar{F}}) \to
\omega_{F_B}$. This  implies
that  the divisor $\alpha_1+\alpha_2+\alpha_3+\alpha_4-h_2$ of $F_B$ is
effective in $F_B$. To obtain the result, we will show that its intersection
with $h_2^2$ is  zero. Remarking from lemma \ref{classFB} that on $F_B$ we have
$h_2^3=24$, we will end the proof after the following lemma:  
\begin{lemma}\label{classalphai}
The image of the class $\alpha_i$ in the Chow ring of $\Gw(2,W)$ is
$h_2.c_2.2(h_2^2-c_2)$, and in the Chow ring of $F_B$ we have $\alpha_i.h_2^2=6$.
\end{lemma}
\proof First, chose an hyperplane $H'$ such that $B=\Gw(3,W)\cap H' \cap H_i$. Let
$F_{H'}$ and $F_{H_i}$  be the variety of lines included respectively in
$\bar{H'}$ and $\bar{H_i}$, and denote by $Y$ the Grassmannian $\Gw(2,p_i^{\bot})$. We
proved in remark \ref{warningFpi} that $\alpha_i$ is reprensented by the
$2$-dimensional part of $Y\cap F_B$ for some point $p_i$ on the
conic $C_i$. From the definition of $B$, we have $F_B=F_{H'}\cap F_{H_i}$. The
variety $F_{H'}$ and $F_{H_i}$ reprensent the class $c_2(\KB)$ in the Chow
ring of $\Gw(2,W)$, but the intersection  $Y\cap F_{H_i}$ has
codimension one in $Y$. Indeed, for any $l$ in $Y$
the point $p_i\wedge l$ is already in $\bar{H}_i$.  The restriction of the
 sheaf $(K_2)^{\bot}$ to $Y$ has a section given by $p_i$,
and the intersection $\Gw(2,p_i^{\bot})\cap F_{H_i}$ is the vanishing locus of
a section of $\frac{K_{2|Y}^{\bot}}{p_i+K_{2|Y}}(h_2)$. Remarking that  $Y$
represent the class $c_2$ in $\Gw(2,W)$, we conclude that the class $\alpha_i$
is equivalent to $h_2.c_2.c_2(\KB)$ in the Chow ring of $\Gw(2,W)$.

Using the remark \ref{chowGw26} we compute that $h^2_2.\alpha_i=6$. This ends the proof
of lemma \ref{classalphai}, and also of lemma \ref{equiv}. \fin

\bigskip
In conclusion, we have the following:
\begin{theorem}\label{FanoB}
 The variety of lines  included in a generic double hyperplane section of $\Gw(3,6)$
  is an hyperplane section of $\PP_1\times\PP_1\times\PP_1\times\PP_1$.
\end{theorem}
\proof In lemma \ref{equiv} we proved that $h_2\sim \sum_{i=1}^4\alpha_i$, but
we have $h^0(\OO_{\bar{F}}(1,1,1,1))=15$ and $h^0(\OO_{\Gw(2,W)}(h_2))=14$, so the morphism
$\psi: F_B \to \bar{F}$ is an isomorphism. \fin
\section{The Chow ring of $B$}
We keep the notations of section \ref{notations}. We will study here the Chow ring of
$B$. The $4$-dimensional variety $B$ is a generic double hyperplane section of
$\Gw(3,6)$. Its Chow ring in codimension $2$ will appear to be surprisingly bigger than
 the codimension $2$ part of the Chow ring of $\Gw(3,6)$. To study it, we will first
 compare it to the Chow ring of $I$ via the projection $p_2: I \to B$.
\begin{lemma}\label{quasifini}
The variety  $B$ contains $16$ quadric cones of dimension $2$. For any point $m$ different
from the $16$ vertex of those cones, the fiber $p_2^{-1}(\{m\})$ has length $4$.
\end{lemma}
\proof Let $m$ be any point of $B$. The  lines containing $m$ and included in
$B$ are the lines containing $m$ and included in the tangent cone of $B$. The tangent cone
to $\Gw(3,W)$ is a cone over a veronese surface. The smoothness of $B$ at $m$ implies that
the tangent cone of $B$ at $m$ is a cone over the intersection of a veronese surface by a
$\PP_3$. So it is  either $4$ lines or a $2$ dimensional quadric cone with
smooth basis.

Now, let $\Gamma$ be such a  cone included in $B$, we can deduce from lemma 3.5
of \cite{I} that there is a point $e$ of $\PP(W)$ included in all the planes
$\pi_u$ for $u$ in $\Gamma$. So $\Gamma$ is the intersection of the
$3$-dimensional quadric 
$Q_e=\{e\wedge l | l\in \Gw(2,\frac{e^{\bot}}{e})\}$ with $2$ hyperplanes containing
$B$. This proves that one of them contains $Q_e$. We deduce from the proposition
3.3 of \cite{I} that this hyperplane must be one of the $H_i$ defined in section
\ref{fibres} and that $e$ is on one of the four conics $C_i$. But the
description of the incidences $Z_i$ in proposition \ref{modelZi} proves that
there are only $4$ cones for each $C_i$. \fin
\begin{description}
\item[\sc Notations:] \ \\
  {\it In the sequel, we will denote the
Chow ring of a variety $X$ by $A_X$. Let $c'_1=h_3$, $c'_2$ and $c'_3$ be the
Chern classes of the tautological   quotient $Q_3$ of $W\otimes
\OO_{\Gw(3,W)}$. Denote by $h'_3$ the class   $p_2^* h_3$ in $A_I$, and let
$a_i$ be the second Chern    class of the bundle $E_i$.}
\end{description}
First recall  from the section \ref{notations} and the lemma \ref{classFB}, that the Chow
ring of the incidence $I$ is:
\begin{lemma}
 The Chow ring of $I$ is
 $$A_I=\frac{A_{F_B}[h'_3]}{(h_3^{'2}-2.h_2.h'_3+\frac{4}{3}.h_2^2)}$$ 
\end{lemma}
So we have to get more informations on $F_B$.
\begin{lemma}\label{blowup}
  For each choice of $i$, the variety $F_B$ can be identified with the blow up of
  $\PP_1\times\PP_1\times \PP_1$ in an elliptic sextic curve, where $v_i$ is the
  exceptionnal divisor. 
\end{lemma}
The choice of one of the $\PP_1$ (ie: the marking of one of the conics $C_i$) in theorem
\ref{FanoB} enable us to consider 
the variety $F_B$ as an incidence:
$$\{(x,h) | x \in \PP_1\times\PP_1\times \PP_1,h\in L, x\in h\}$$
where $L$ is a marked subspace of
$H^0(\OO_{\PP_1\times\PP_1\times \PP_1}(1,1,1))$. In other words, for each $i$,
the variety $F_B$ with le line bundle $h_2-\alpha_i$ is identified with the blow
up of $\PP_1\times\PP_1\times \PP_1$ in an elliptic sextic curve. By
construction $|\alpha_i|$ is the marked pencil of hyperplanes.

To show that $v_i$ is the exceptionnal divisor, we need to prove that $2\alpha_i+v_i=
h_2$. Take any subset $\{i,j,k\}$ of $\{1,2,3,4\}$ of cardinal $3$, and chose  a point
$p_j$ on $C_j$ and another one $p_k$ on $C_k$ such that $\PP(p_j^{\bot}\cap
p_k^{\bot})\cap C_i=\emptyset$. The divisor $v_i$ was constructed in proposition
\ref{alphai} as the set of lines included in the hyperplane section $\bar{H}_{u_i}$ of
$B$. So $\delta \in v_i$ implies that $\Delta \subset  \bar{H}_{u_i}$ which implies that
the line $\PP(L_\delta)$ intersect $C_i$, so it can't be in $\PP(p_j^{\bot}\cap
p_k^{\bot})$, and we have $\alpha_j.\alpha_k.v_i=0$. So set theoretically, $v_i$ is the
exceptionnal divisor, but we proved in lemma \ref{reduit} that $v_i$ is
reduced.  \fin
\begin{lemma}\label{EiI}
  The projection $p_{1*}(p_2^*a_i)$ of the second Chern class of $p_2^*E_i$ to
  $A_{F_B}$ is   $h_2-\alpha_i$. Furthermore, we have in $A_I$ the equality:
  $p_2^*a_i=h'_3.(h_2-\alpha_i)+h_2(2\alpha_i-h_2)$   
\end{lemma}
\proof We know from the lemma \ref{blowup} and proposition \ref{alphai}, that
the first Chern class of $p_{1*}(p_2^* E_i )$ is $-\alpha_i-v_i=\alpha_i-h_2$. So
the first assertion is a direct consequence of Riemman-Roch-Grothendieck's
theorem because $(p_2^* E_i )$ have no higher direct images by $p_1$. Indeed, we
can compute in $A_I$ using  the equality: $\omega_{p_1}=\OO_I(2h_2-2h_3)$
obtained in  {\S}\ref{notations}, and  it gives that the first Chern class of
$p_{1!}(p_2^* E_i )$ is $-p_{1*}p_2^*a_i$.  

The second assertion is obtained by the evaluation map: $p_1^*p_{1*} p_2^* E_i
\to p_2^* E_i$. We can compute its cokernel as we done in lemma \ref{Brel}, but
the vanishing of  $R^1p_{1*}(p_2^*E_i)$ gives the following exact sequence:
$$0 \to p_1^*p_{1*} p_2^* E_i \to p_2^* E_i \to \OO_I(-h_3')\otimes
p_1^*R^1p_{1*}(p_2^*(E_i(-h_3))) \to 0$$
so we have by relative duality the extention:
\begin{equation}
  \label{p2Ei}
0 \to \OO_I(p_1^*\alpha_i-p_1^*h_2) \to p_2^* E_i \to
\OO_I(p_1^*h_2-p_1^*\alpha_i-h'_3) \to 0   
\end{equation}
which gives the computation of $p_2^*a_i$. \fin
\begin{lemma}\label{pullback}
  Let $\gamma\in A^2_B$, then the class $p_2^* \gamma$ can be written in $A_I$
  by $h'_3.\gamma_0+\gamma_1$ with   $\gamma_i$ in $A^{i+1}_{F_B}$, and where
  $\gamma_i$ is in the vector space generated by
  $h_2^i\alpha_1,\dots,h_2^i\alpha_4$. More precisely, we have:
  $2h_2.\gamma_0+\gamma_1\in\ratio .h_2^2$   
\end{lemma}
\proof We can first find classes $\gamma_i$ in $A^{i+1}_{F_B}$ such that
$p_2^*\gamma=h'_3.\gamma_0+\gamma_1$.  Now remark that  $A^3_B$ is one
dimensional, so $h_3.\gamma$ is propotional to 
$h^3_3$. But from the relation $h_3^{'2}=2.h_2.h'_3-\frac{4}{3}.h_2^2$, the
class $p_{1*}h_3^{'3}$ is proportional to $h_2^2$. So the class
$p_{1*}p_2^*\gamma$ is also in $\ratio.h_2^2$, and we have
$2h_2.\gamma_0+\gamma_1\in\ratio .h_2^2$. We can now conclude with the
description of $F_B$ in lemma \ref{blowup} that $\gamma_1$ and $h_2$ are in the
vector space generated by $\alpha_1,\dots,\alpha_4$, which gives the lemma. \fin
\begin{lemma}\label{gensA2B}
  The classes $(a_1,a_2,a_3,a_4)$ form a basis of the vector space  $A^2_B$. We have in
  $A_B$ the   relation $2(a_1+a_2+a_3+a_4)=3h_3^2$. 
\end{lemma}
\proof From the lemma \ref{quasifini}, the map  $p_2^*: A_B^2 \to A^2_I$ is
injective, so  from the lemma \ref{pullback}, $A_B^2$ is generated by the
set of classes $\{p_{2*}(h'_3.\alpha_i),p_{2*}(h_2\alpha_i)\}_{i=1\dots 4}$.
 As the picard group of $B$ is generated by $h_3$, all the classes
$p_{2*}(h'_3.\gamma_0)$  are proportional to $h_3^2$. Now, we
use the relation $\frac{4}{3}.h_2^2=2.h_2.h'_3-h_3^{'2}$ to eliminate $h_2^2$ in the
expression of $p_2^*a_i$ found in lemma \ref{EiI}. So, we obtain that
$a_i$ is in the affine space $\frac{1}{2}.p_{2*}(h_2\alpha_i)+\ratio.h_3 ^2$.
So we have
proved that $A^2_B$ is generated by $h_3^2,a_1,a_2,a_3,a_4$. Furthermore, as we have
$p_{1*}p^*_2a_i=h_2-\alpha_i$, which is a free family in $A^1_{F_B }$, the family
$(a_1,\dots,a_4)$ is free in $A^2_B$.

To obtain the relation with $h_3^2$, we substitute the expression of $p_2^*a_i$ of lemma
\ref{EiI} in the relation found in lemma \ref{equiv}. We eliminate 
$h'_3.h_2$ with the 
relation $\frac{4}{3}.h_2^2+h_3^{'2}=2.h_2.h'_3$, and we obtain: $\sum\limits_{i=1}^4 p_2^* a_i =
\frac{3}{2} h_3^{'2}$. \fin

So we are now ready to describe the Chow ring of $B$. 
\begin{proposition}\label{ChowB}
  The Chow ring of $B$ is $\ratio[h_3,a_1,a_2,a_3,a_4]/\II$
  where $\II$ is generated by $3h_3^2-2.\sum_{i=1}^4a_i$,
  $(8h_3.a_i-3h_3^3)_{i\in\{1,\dots,4\}}$, $(8.a_i.a_j-h_3^4)_{i\neq
    j,(i,j)\in\{1,\dots,4\}^2}$.  (The class of a point is $\frac{a_1.a_2}{2}$,
  the class of the veronese $V_i$ is $2a_i-\frac{1}{2}h_3^2$, and $[V_i]^2$ is
  $4$ points).   
\end{proposition}
\proof As $A^i_B$ is known to be one dimensional for $i\in\{0,1,3,4\}$, we just need to compute the
relations, and it can be done by calculating the degrees, because we have found the
structure of $A^2_B$ in lemma \ref{gensA2B}. The relations $(8h_3.a_i-3h_3^3)$
are consequences of the degree of $\Gw(3,6)$ (ie 16) and the fact that $a_i$ can
be represented by the union of a quadric $Z_{i,p_i}$ and the veronese surface $V_i$.

Now remark that the class of  $V_i$ is in the vector space generated by $a_i$ and
$h_3^2$. Indeed, we have $p_2^*[V_i]=h'_3.v_i+\gamma_1$ and the lemmas
\ref{blowup} and \ref{pullback} give $v_i=h_2-2\alpha_i$ and $\gamma_1\in 
h_2^2.\ratio+h_2\alpha_i.\ratio$, so $[V_i]$ is in the
vector space generated by $h_3^2$ and $a_i$, and then $[Z_{i,p_i}]$
also. Computing their degree, we have: $[V_i]=2a_i-\frac{1}{2}h_3^2$ and
$[Z_{i,p_i}]=\frac{1}{2}h_3^2-a_i$. The last relations are then consequences of
the fact that we have: $[V_i].[V_j]=0$ for $i\neq j$. \fin

NB:

\bigskip
It could be usefull to have the link with the ring of $\Gw(3,6)$ and the Chern
classes of $Q_3$, so we state the following:
\begin{remark}\label{chowGw36}
  Still denote by $c'_i$ the Chern classes of the tautological quotient $Q_3$. The Chow
  ring of $\Gw(3,6)$ is 
  $$\ratio[c'_1,c'_2,c'_3]/(c_3^{'2},c_2^{'2}-2c'_1c'_3,c_1^{'2}-2c'_2)$$
  On $B$ we have the additional relation $c'_1c'_2=4c'_3$. In particular the
  rank of $A^2_{\Gw(3,6)}$ is only one.
\end{remark}
\proof The Chow ring of $\Gw(3,6)$ is given by
$$\ratio[x,y,z]/((xyz)^2,x^2y^2+y^2z^2+z^2x^2,x^2+y^2+z^2)$$
where $c'_1=x+y+z$, $c'_2=xy+yz+zx$ and $c'_3=xyz$. \fin

\section{Application to quadratic normality}
In this part, we  explain the link with the congruence of lines found in
\cite{MdP}. They started with the intersection $\Gamma$ of $G(2,6)$ by a very particular
$\PP_{11}$ (NB: it was proved in \cite{MM}, that the choice of this $\PP_{11}$ is unique
up to the action of $GL_6$). Then they chosed a general quadric containing a fixed
subscheme of $\Gamma$ (which contains the singular locus of $\Gamma$) to obtain a
reducible intersection with $\Gamma$. Then they checked with Macaulay2 that one of these
irreducible components is smooth of degree $16$ and sectional genus $9$.

Here we prove that the generic Fano $4$-fold of genus $9$ can be obtained by their
construction, and that the choice we need to do is generically finite. More
precisely, the choice of $\Gamma$ correspond to the choice of 
a tangent hyperplane $H$ to $\Gw(3,6)$, and the choice of the quadric will
correspond to the choice of a non zero element in $|\OO_{\bar{H}}(h_3)|$.
\begin{proposition}\label{congruence}
  For any choice of $i$, the bundle $E_i(h_3)$ gives an embedding of the Fano manifold $B$
  in the grassmannian $G(2,6)$ as the congruence of lines constructed in
  \cite{MdP} (Theorems 8,9,10).  
\end{proposition}
\proof Their description of the congruence is in terms of equations, so we need
to make an adapted choice of coordinates to get the link. Choose a basis
$w_0,\dots,w_5$ of $W$ such 
that $\omega\dual=w_0\wedge w_3+w_1\wedge w_4+w_2\wedge w_5$. Let $A$ and $B$ be the
vector spaces generated respectively by $w_0,w_1,w_2$ and $w_3,w_4,w_5$, in particular,
the form $\omega$ gives an identification between $A\dual$ and $B$. The decomposition
$W=A\oplus B$ gives a decomposition of $\bigwedge\limits^3 W$. So we can represent an
element of $\bigwedge\limits^3 W$ as in \cite{I} by $(a,X,Y,b)$, with $a\in \wedge^3 A$,
$b \in \wedge^3 B$, $X\in Hom(A,B)$, $Y\in Hom(B,A)$.  The equations of $G(3,6)$ are
 as: $\wedge^2 X=aY$, $\wedge^2 Y=bX$, $YX=ab.I_3$,
and to obtain $\Gw(3,6)$ we had the linear relations $X= \T X$ and $Y= \T Y$, so we take $X=\left(
  \begin{array}{ccc}
    x_0&x_1&x_2\\
    x_1&x_5&x_3\\
    x_2&x_3&x_4
  \end{array}\right)$ and $Y=\left(
  \begin{array}{ccc}
    y_0&y_1&y_2\\
    y_1&y_5&y_3\\
    y_2&y_3&y_4
  \end{array}\right)$. (Cf \cite{I}).

Now we need to choose an hyperplane $H$ as in section \ref{sectionhyp}. So we
consider the hyperplane $y_{5}=y_{2}$. It is  tangent to $\Gw(3,6)$ in
$w_0\wedge w_1\wedge w_2$ and contains $w_3\wedge w_4\wedge w_5$. Furthermore,
the conic decribed in lemma \ref{LC} is  parametrized by
$\lambda^2w_0+\lambda\mu w_1+\mu^2w_2$,
so  the incidence $Z_H$ is given in $\PP_1\times \bar{H}$ by the equations:
$$a=0, X.{\tiny\left(\begin{array}{c}
  \lambda^2\\
  \lambda\mu\\
  \mu^2
\end{array}\right)}={\tiny\left(\begin{array}{c}
  0\\
  0\\
  0
\end{array}\right)},
{\footnotesize\left(
\begin{array}{ccc}
  -\mu & \lambda & 0\\
  0 & -\mu &\lambda
\end{array}
\right)}.Y={\footnotesize \left(
  \begin{array}{ccc}
    0&0&0\\
    0&0&0
  \end{array}\right)
} $$
The 6 dimensional vector space $H^0(\II_{Z_H}(1,1))$ is generated by
$$(-a.\lambda,-a.\mu,\lambda y_1-\mu y_0, \lambda y_2-\mu y_1, \lambda y_3 - \mu
y_2, \lambda y_4 -\mu y_3)$$
So, the sheaf $\EE$ constructed in
proposition \ref{defE} is the image of the
map:
$$
6\OO_{\bar{H}}(-1) \begin{array}[b]{c}
{\left(
\begin{array}{cccccc}
  a & 0 & y_0 & y_1 & y_2 & y_3\\
  0 & a & -y_1 & -y_2 & -y_3 & -y_4
\end{array}\right)}\\
\overrightarrow{\hspace{14em}}
\end{array} 2 \OO_{\bar{H}}
$$
which is exactly the map of \cite{MdP} th 9.
Now,  a generic hyperplane section of $\Gw(3,6)$ gives a linear relation
between $a,(x_i),(y_i),b$ which is exactly the relation 29 of \cite{MdP}, which
gives the identification with their congruence of lines. \fin

So we can remark that the computations they made in affine charts of $G(2,6)$ can be
globally done in $\Gw(3,6)$.

\subsection{Geometry of the focal locus of Mezzeti-de Poi's congruence}
Here we chose one of the $4$ bundles, say $E_1$. Denote by $h$ the bundle
$\OO_{Proj(E_1(h_3))}(1)$. We proved in lemma \ref{acyEi} that the linear system
$|h|$ gives a morphism from $Proj(E_1(h_3))$ to $\PP_5$.
\begin{description}
\item{\sc Notations:} \\
  {\it Denote by $r$ the above projection from $Proj(E_1(h_3))$ to $\PP_5$, and by $p_B$ the
    projection from $Proj(E_1(h_3))$ to its basis $B$. }
\end{description}
Following \cite{MdP},
the focal locus of Mezzeti-de Poi's congruence is  defined like this:
\begin{lemma}\label{classR}
  The morphism from $r: Proj(E_1(h_3))\to\PP_5$ is birational. The class of the
  exceptional divisor ${\cal R}$ in $ Proj(E_1(h_3))$ is $4h-h_3$.
\end{lemma}
\proof The degree of $r$ is given by the length of the degeneracy locus of a
generic map $5\OO_B \to E_1(h_3)$, so it is the fourth segre class of $E_1$,
which is the class of one point from proposition \ref{ChowB}. So $r$ is birational,
and the class of the exceptional divisor follows. \fin

\bigskip
The manifold $B$ gives  the
family $LB=\{\PP(E_{1,u})| u\in B\}$ of lines in $\PP_5$. From the lemma \ref{classR},
those lines are quadrisecant to $X$. 
\begin{description}
\item{\sc Notations:} \\
{\it Denote by $X$ the focal locus of $B$, (ie $X=r({\cal R}$)). Recall from \cite{MdP}
  that $X$ has dimension $3$, degree $6$, and is 
singular along a rational smooth cubic curve $C$.} 

\end{description}
 From the lemma
\ref{classR}, any line of this family  intersects $X$ in
length $4$ or is included in $X$. We can now describe easily the normalisation
$\tilde{X}$ of $X$, and  
note that we have as in the Palatini case some kind of duality between $X$ and a
family of plane 
cubics included in $X$, but here it breaks over the singular locus of $X$:
\begin{proposition}\label{descX}
  The focal locus $X$ of Mezzeti-de Poi's congruence is a projection of
  $\PP_1\times\PP_1\times\PP_1$ from a  line. The variety of pencils of
   lines belonging to $LB$ is $\PP_1\times\PP_1\times\PP_1$ blown up in the
  elliptic curve of degree $6$ which is the double cover of the cubic curve $C$.
\end{proposition}
\proof We pullback the situation to the incidence  point/lines of $B$ via the
projection $p_2: I \to B$. We still denote by $p_2$ the projection from
$Proj(p_2^*(E_1(h_3)))$  to $Proj(E_1(h_3))$. In 
the proof of lemma  \ref{EiI}, we noticed  the extension (\ref{p2Ei}), which gives a
surjection: 
$$p_2^*(E_1(h_3)) \to\hspace{-0.8em}\to \OO_I(p_1^*h_2-p_1^*\alpha_1)$$
This surjection gives an embedding of $I$ in $Proj(p_2^*(E_1(h_3)))$.  The
restriction of $p_2^*(h)$ to  $I$ is $p_1^*(h_2-\alpha_1)$, but it is also
 $p_1^*(\alpha_2+\alpha_3+\alpha_4)$ by lemma \ref{equiv}. So the linear system
 $|p_2^*(h)|$ contracts the fibers of $p_1:I\to F_B$, and the image of $I$
 coincide with    the  image of $F_B$ by $\alpha_2+\alpha_3+\alpha_4$. In
 conclusion, the map $r$ contracts the divisor $p_2(I)$ to the projection of
 $\PP_1\times \PP_1\times \PP_1$ from a  line. \fin
 
\bigskip 
We can also obtain a description of the plane curves in $X$ related to this family of
lines:
\begin{remark}\label{12dtes}
Any point of $F_B$ gives a plane in $\PP_5$ intersecting $X$ in a point and a plane
cubic. Those plane cubics are all singular in a point of $C$. Only
 $12$ lines included in $X$ are in $LB$.
\end{remark}
\proof A point  $\delta\in F_B$ correspond to a line $\Delta$ included in $B$. We already
proved  in {\S}2 that the restriction of $E_1$ to $\Delta$ is always $\OO_\Delta \oplus
\OO_\Delta(-1)$, so the lemma \ref{classR} implies that the intersection of ${\cal R}$
with $p_B^{-1}(\Delta)$ is given by a section of $(S_4(E_1))(3h_3)$. Hence this
intersection contains the exceptionnal divisor of $\PP(\OO_\Delta\oplus
\OO_\Delta(1))$, so 
the plane $r(p_B^{-1}(\Delta))$ intersect $X$ in a cubic curve and a unique
residual point included in all the lines $r(\PP(E_{1,p}))$ for any point $p$ of
$\Delta$. Note  that the exceptional  divisor ${\cal R}$ is the image of the composition:
\begin{equation}
  \label{factorisation}
 I\subset Proj(p_2^*(E_1(h_3))) \to Proj(E_1(h_3))  
\end{equation}
so from lemma \ref{quasifini} the
 projection ${\cal R} \to B$ is finite of degree $4$ except over $16$ points of  $B$, but
 we will notice later that $4$ of them are contracted to points on $X$.
To understand why the  plane cubic described above is singular, we first notice that it is
the image of the following curve $T_\Delta$ in $F_B$: The closure in $F_B$ of the lines included in
$B$ intersecting $\Delta$ and different from $\Delta$. A general $\Delta$, intersect the
hyperplane section $\bar{H}_{u_1}$ in a single point $b$ which is not on the veronese
$V_1$. From the corollaries \ref{pinceau} and \ref{EiZi}, there is a single quadric
$Z_{1,p_1}$ containing $b$. The $2$ lines containing 
$b$ and included in this quadric are included in $\bar{H}_{u_1}$. So, from remark
\ref{secVi}, they correspond in $F_B$ to points on the exceptional divisor $v_1\subset F_B$. 
In conclusion, for a general $\Delta$, the curve $T_\Delta$ intersects the exceptional divisor
$v_1$ in $2$ points. The image of those $2$ points of $F_B$ must be a single point of $C$
because they are in the line $r(\PP(p_B^{-1}(b)))$, and this line intersects $C$ in at
most one point because $b$ is not in the veronese $V_1$. So the plane cubic is singular at
this point of $C$.

For the same reason, if $b$ is one of the $4$ vertex of the cones in  $\bar{H}_{u_1}\cap
B$, $r(\PP(p_B^{-1}(b)))$ is contracted to a point on the curve $C$, that's why only $12$
lines of $X$ belong to the family $LB$. \fin

\bigskip
As  all lines in those planes are trisecant to $X$,  we  have from the above remark that: 
\begin{remark}
The triple locus of the projection of $\PP_1\times\PP_1\times \PP_1$ from a
generic plane is reducible, because it contains the union of $4$ lines with a common
point.     
\end{remark}
\bigskip
We give now some details on the marked ``virtual'' section related to the anormality of
$X$.
\begin{remark}
  The ramification of the morphism $p_2: I \to B$ is $p_1^{-1}(\Sigma)$ for some surface
  $\Sigma$ in $F_B$. The canonical sheaf of $\Sigma$ is $\omega_\Sigma=\OO_\Sigma$, and
  $\Sigma$ contains the $16$ rational curves parameterizing  the $16$ cones of $B$. The
  image of $\Sigma$ into $X$ is a section of $\OO_X(2)$ that is not in a quadric of
  $\PP_5$. 
\end{remark}
\proof The canonical divisor of $I$ can be computed from {\S}2, and it gives that  ${\cal R}
\sim p_1^* h_2$. The $16$ contracted curves must be in ${\cal R}$, so we have only to
prove that the image of $\Sigma$ in $X$ is not in a quadric of $\PP_5$. But the answer is
general for those anormality constuctions: As ${\cal R}$ is a divisor in the projective
bundle $Proj(E_1(h_3))$, we have the exact sequence:
$$0 \to \omega_{p_B} \to \omega_{p_B}({\cal R}) \to \omega_{\cal R} \otimes \omega_B\dual
\to 0$$
where $\omega_{p_B}$ is the relative dualising sheaf of $p_B: Proj(E_1(h_3))\to B$. The
ramification of the restriction of $p_B$ to ${\cal R}$ is the zero locus of a section of $
\omega_{\cal R} \otimes \omega_B\dual$ which gives a non zero element of
$H^1(\omega_{p_B})=\complex$. But we saw in the proof of remark \ref{12dtes} that the map $I\to
B$ factors through ${\cal R}$ (Cf the composition (\ref{factorisation})), so we have the
result because the above sequence is also:
$$0 \to \omega_{p_B} \to \OO_{Proj(E_1(h_3))}(2h) \to \OO_{\cal R}(2h) \to 0 $$
\fin
\subsection{The rank $2$ reflexive sheaf on $\PP_5$}
If we take a general double hyperplane section of $X$, we obtain a smooth elliptic curve
of degree $6$. So from Serre's construction, it is the vanishing locus of a section of a
rank $2$ vector bundle on $\PP_3$, unique up to isomorphism. Curiously, there  a way to
globalise this construction over $\PP_5$. In this part we will construct an
$SL_2$-equivariant rank $2$ reflexive sheaf on $\PP_5$ using classical
techniques developped in mathematical instanton studies. (Cf \cite{Ba}, \cite{Tj})
Still consider  vector spaces $L$  and $V$ of respective dimension $2$ and $6$.
 This construction will be essentially unique up to the $SL_6$ action. Indeed,
 it could be constructed from a tangent hyperplane to $\Gw(3,6)$ or like this:

Let's recall from \cite{MM} that $S_2 L\otimes \bigwedge\limits^2 V$ have an $SL_2\times
SL_6$ orbit  made of
net of alternating forms of constant rank $4$.  So let's $\beta \in S_2 L\otimes
\bigwedge\limits^2 V$ be the  element of this orbit. This element was considered
in \cite{MdP} to construct a $\PP_{11}$ containing their congruence, for
instance, take the following:
{\small$$\left(
\begin{array}{cccccc}
  0&u^2&2uv&v^2&0&0\\
  -u^2&0&0&0&0&0\\
  -2uv&0&0&0&0&u^2\\
  -v^2&0&0&0&0&2uv\\
  0&0&0&0&0&v^2\\
  0&0&-u^2&-2uv&-v^2&0
\end{array}
\right)$$}
We can remark that $\beta$ viewed as an element of $\bigwedge\limits^2 (L\otimes V)$ has
rank $6$ because it can be represented by:
{\small$$\left(
\begin{array}{cccccccccccc}
  0&0&1&0&0&1&0&0&0&0&0&0
  \\0&0&0&0&1&0&0&1&0&0&0&0
  \\-1&0&0&0&0&0&0&0&0&0&0&0
  \\0&0&0&0&0&0&0&0&0&0&0&0
  \\0&-1&0&0&0&0&0&0&0&0&1&0
  \\-1&0&0&0&0&0&0&0&0&0&0&0
  \\0&0&0&0&0&0&0&0&0&0&0&1
  \\0&-1&0&0&0&0&0&0&0&0&1&0
  \\0&0&0&0&0&0&0&0&0&0&0&0
  \\0&0&0&0&0&0&0&0&0&0&0&1
  \\0&0&0&0&-1&0&0&-1&0&0&0&0
  \\0&0&0&0&0&0&-1&0&0&-1&0&0
\end{array}\right)
$$}
Denote by $W$ the six dimensional image of this map, we will identify $W$ with
its dual via the induced alternating form. The inclusion $W\subset L\otimes V$
gives a map $\beta'$ from $W$ to $L\otimes \OO_{\PP_5}(1)$, and  we can construct a complex:
\begin{equation}
  \label{quasimonade}
L\otimes \OO_{\PP_5}(-1) \stackrel{^t\beta'}\longrightarrow  W\otimes \OO_{\PP_5}
\stackrel{\beta'}{\longrightarrow} L\otimes \OO_{\PP_5}(1)   
\end{equation}
exact on the left, with middle cohomology a rank $2$ reflexive sheaf $K$, and
with right  cohomology a sheaf $\LL$ supported on a smooth rational cubic curve
$C$ isomorphic to 
$\OO_{\PP_1}(4)$. So we can compute from the complex (\ref{quasimonade}) some
of its invariants: 
\begin{corollary}
  The sheaf $K$ has rank $2$, $c_1K=0$, $c_2K=2$, $c_3K=0$, $c_4K=-15$. Its
  singular locus is the cubic curve $C$, and we have
  $H^0K(1)=0$ and $H^0K(2)=13$, and $H^1K=1$.
\end{corollary}
we have now a way to structurate the variation of the focal locus with respect
to the choice of the quadric in th8 of \cite{MdP}. They are the vanishing locus
of sections of the same bundle: $K(2)$.


\begin{thebibliography}{}
\bibitem[Ba]{Ba}  Barth, W: Irreductibility of the space of
mathematical instanton bundles with Rank 2, $c_2=4$. Math.Ann {\bf 258}, 81-106 (1981)
\bibitem[G-P]{GP} L.Gruson - C.Peskine: Courbes de l'espace
projectif,vari\'et\'es de s\'ecantes. Enumerative geometry and classical
algebraic geometry. Nice (1981) Prog Math {\bf 24} Birkh\"auser.Boston
1982
\bibitem[I]{I} A. Iliev: The $SP_3$-grassmannian and duality for prime Fano threefolds of
  genus $9$. Manuscripta math. {\bf 112}, 29-53 (2003)
\bibitem[I-R]{IR} A. Iliev - K. Ranestad: Geometry of the Lagrangian
  Grassmannian LG(3,6) with applications to Brill-Noether
  loci. Mich. Math. Journal {\bf 53}, 383-417 (2005)
\bibitem[K]{K} A. Kuznetsov: Hyperplane sections and derived
  categories. Izvestiya Mathematics {\bf 70}:3, 447-447 (2006) 
\bibitem[M]{M} L.Manivel: Configuration of lines and models of Lie
  algebras. Journal of Algebra {\bf 304}, 457-486 (2006)
\bibitem[M-M]{MM} L. Manivel - E. Mezzeti: On linear spaces of skew-symmetric matrices of
  constant rank. Manuscripta math. {\bf 117}, 319-331 (2005)
\bibitem[M-dP]{MdP} Mezzeti - de Poi: Congruences of lines in $\PP_5$, quadratic
  normality, and completely exceptional Monge-Amp\`ere equations. Geom Dedicata {\bf 131}, 213-230 (2008)
\bibitem[Tj]{Tj}  Tjurin, A.N: On the superpositions of mathematical
instantons. In Artin, Tate, J.(eds) Arithmetic and geometry. Prog.Math
{\bf 36}, 433-450 (1983) . Birkh\"auser.
\bibitem[W]{W} J. Weyman: Cohomology of vector bundles and syzygies. {\bf 149} (2003)
  Cambridge   Tracts in Mathematics.
\end{thebibliography}
\end{document}